\documentclass[10pt,twoside]{amsart}
\advance\oddsidemargin by -1.0cm
\advance\evensidemargin by -1.0cm
\advance\topmargin by -1.0cm 
\usepackage{amssymb}

\theoremstyle{plain}

\newtheorem*{lem}{Lemma}
\newtheorem{thm}{Theorem}
\newtheorem{prop}{Proposition}
\theoremstyle{remark}

\newtheorem*{rem}{Remark}
\newtheorem{exa}{Example}
\newcommand{\bdm}{\begin{displaymath}}
\newcommand{\edm}{\end{displaymath}}
\newcommand{\ba}[1]{\begin{array}{#1}}
\newcommand{\ea}{\end{array}}
\newcommand{\bea}[1][]{\begin{eqnarray#1}}
\newcommand{\eea}[1][]{\end{eqnarray#1}}

\newcommand{\bb}{\begin{bundle}}
\newcommand{\eb}{\end{bundle}}
\newcommand{\bcen}{\begin{center}}
\newcommand{\ecen}{\end{center}}
\newcommand{\btab}{\begin{tabular}}
\newcommand{\etab}{\end{tabular}}

\newcommand{\x}{\times}

\newcommand{\ox}{\otimes}

\newcommand{\ra}{\rightarrow}
\newcommand{\lra}{\longrightarrow}

\newcommand{\lan}{\left\langle}
\newcommand{\ran}{\right\rangle}

\newcommand{\de}{\ensuremath{\mathrm{d}}}
\newcommand{\vol}{\ensuremath{\mathrm{vol}}}
\newcommand{\restr}{\arrowvert}
\newcommand{\C}{\ensuremath{\mathbf{C}}}
\newcommand{\R}{\ensuremath{\mathbf{R}}}

\newcommand{\N}{\ensuremath{\mathbf{N}}}

\newcommand{\eps}{\ensuremath{\varepsilon}}

\newcommand{\Conto}[1]{\ensuremath{\mathcal{C}^0(#1)}} 
\newcommand{\Contoinfty}[1]{\ensuremath{\mathcal{C}^0_{\infty}(#1)}} 
\newcommand{\Pol}[1][V]{\ensuremath{\mathcal{P}(#1)}}   
\newcommand{\Aff}[3][]{\ensuremath{#2^{#1}[#3]}}        
\begin{document}
\setcounter{equation}{0}
%
%
\thispagestyle{empty}
%
%
\date{\today}
\title{The Gaussian measure on algebraic varieties}
%
%
%
\author{Ilka Agricola}\author{Thomas Friedrich}
\address{\hspace{-5mm} 
{\normalfont\ttfamily agricola@mathematik.hu-berlin.de}\newline
{\normalfont\ttfamily friedric@mathematik.hu-berlin.de}\vspace{2mm}\newline
Institut f\"ur Reine Mathematik \newline
Humboldt-Universit\"at zu Berlin\newline
Sitz: Ziegelstr. 13 A\newline
D-10099 Berlin\\
Germany}
\thanks{This work was supported by the SFB 288 "Differential
geometry and quantum physics".}
\begin{abstract}
We prove that the ring $\Aff{\R}{M}$
of all polynomials defined on a real algebraic variety $M\subset\R^n$
is dense in the Hilbert space $L^2(M,e^{-|x|^2}\de\mu)$, where
$\de\mu$ denotes the volume form of $M$ and $\de\nu=e^{-|x|^2}\de\mu$ the 
Gaussian measure on $M$.
\end{abstract}
\maketitle
%
\pagestyle{headings}
%
%
%
\section{Introduction}
The aim of the present note is to prove that the ring $\Aff{\R}{M}$
of all polynomials defined on a real algebraic variety $M\subset\R^n$
is dense in the Hilbert space $L^2(M,e^{-|x|^2}\de\mu)$, where
$\de\mu$ denotes the volume form of $M$ and $\de\nu=e^{-|x|^2}\de\mu$ is the 
Gaussian measure on $M$. In case $M=\R^n$, the result is well known since the
Hermite polynomials constitute a complete orthonormal basis of
$L^2(\R^n,e^{-|x|^2}\de\mu)$.

\section{The volume growth of an algebraic variety and some consequences}
We consider a smooth algebraic variety $M\subset\R^n$ of dimension $d$.
Then $M$ has polynomial volume growth: there exists a constant $C$ depending
only on the degrees of the polynomials defining $M$ such that for any euclidian
ball $B_r$ with center $0\in\R^n$ and radius $r>0$ the inequality
 \bdm
 \vol_d (M\cap B_r)\ \leq\ C\cdot r^d
 \edm
holds (see \cite{Broecker95}). Via Crofton formulas the mentioned inequality
is a consequence of Milnor's results concerning the Betti numbers of an
algebraic variety (see \cite{Milnor64}, \cite{Milnor94}, in which the
stated inequality is already implicitly contained). This estimate yields
first of all that the restrictions on $M$ of the polynomials on $\R^n$ are
square-integrable with respect to the Gaussian measure on $M$.

\begin{prop}
Let $M$ be a smooth submanifold of the euclidian space $\R^n$. Suppose
that $M$ has polynomial volume growth,i.e., there exist constants
$C$ and $l\in\N$ such that for any ball $B_r$
 \bdm
 \vol_d (M\cap B_r)\ \leq\ C\cdot r^l
 \edm
holds. Denote by $\de\mu$ the volume form of $M$. Then:
\begin{enumerate}
\item The ring $\Aff{\R}{M}$ of all polynomials on $M$ is contained in 
the Hilbert space $L^2(M,e^{-|x|^2}\de\mu)$;
\item all functions  $e^{\alpha |x|^2}$ for $\alpha < 1/2$ belong to
$L^2(M,e^{-|x|^2}\de\mu)$.
\end{enumerate}
\end{prop}

\proof 
Throughout this article, denote the distance of the point
$x\in\R^n$ to the origin by $r^2=|x|^2$.
We shall prove that the integrals 
 \bdm
 I_m(M)\ :=\ \int_M r^me^{-r^2}\de\mu < \infty,\quad m=1,2,\ldots
 \edm
are finite. However,
 \bdm
 I_m(M)\ =\ \sum_{j=0}^{\infty} 
 \int_{M\cap (B_{j+1}-B_j)} r^m e^{-r^2}\de\mu 
 \edm
and consequently we can estimate $I_m(M)$ as follows:
 \bdm
 I_m(M)\ \leq\  \sum_{j=0}^{\infty}
 (j+1)^m e^{-j^2}\left[ \vol (M\cap B_{j+1}) - \vol (M\cap B_j) \right]
 \ \leq\  \sum_{r=0}^{\infty}
 (r+1)^m e^{-r^2} \vol (M\cap B_{r+1})
 \, .
 \edm
Using the assumption on the volume growth of $M$ we immediately obtain
 \bdm
 I_m(M)\ \leq\  C\cdot\sum_{r=0}^{\infty}
 (r+1)^{m+l} e^{-r^2}\, .
 \edm
Denoting the summands of the latter  series by $a_r$, we readily see that it
converges, since 
 \bdm
 \frac{a_{r+1}}{a_r}\ =\ \frac{(r+1)^{m+l}e^{-r^2-2r-1}}{(r)^{m+l}e^{-r^2}}\
 =\ \left( \frac{r+1}{r}\right)^{m+l} \frac{1}{e^{2r+1}}\lra 0.
 \edm
A similar calculation yields the result for the functions $e^{\alpha r^2}$
with $\alpha < 1/2$.
\qed

\section{A dense subspace in $\Contoinfty{S^n}$}
The aim of this section is to verify that a certain linear subspace
of $\Conto{S^n}$ is dense therein. Since the family of functions we have in
mind cannot be made into an algebra, we have to replace the standard
Stone-Weierstra\ss{}
argument by something different. The idea for overcoming this 
problem is to use a combination of the well-known
theorems of Hahn-Banach, Riesz and Bochner.

To begin with, we uniformly approximate the function
$e^{-r^2}e^{i\lan  k,x\ran}$ for a fixed vector $k\in\R^n$.
\begin{lem}
Denote by $p_m(x)$ the polynomial
 \bdm
 p_m (x) = \sum_{\alpha=0}^{m-1}i^{\alpha}\lan k,x \ran ^{\alpha}/ \alpha !\, .
 \edm
Then the sequence $e^{-r^2}p_m(x)$ converges uniformly  to
$e^{-r^2}e^{i\lan  k,x\ran}$ on $\R^n$.
\end{lem}
\proof
The inequality
 \bdm
 |\, p_m(x) - e^{i\lan  k,x\ran}\, |\ \leq\ \frac{\|k\|^m \|x\|^m}{m!}\
 e^{\|k\|\cdot\|x\|} 
 \edm
implies (set  $y=\|k\|\cdot \|x\|$)
 \bdm
 \sup_{x\in\R^n} |\,e^{-r^2} p_m(x) -e^{-r^2} e^{i\lan  k,x\ran}\, |\ \leq\ 
 \sup_{0\leq y} \ \frac{y^m}{m!} \ e^{y- y^2/\|k\|^2} \ =: \ C_m \, .
 \edm
Therefore, we have to check that for any fixed vector $k\in\R^n$ the
sequence $C_m$ tends to zero as $m\ra \infty$. For simplicity, denote by
$k$ the length of the vector $k\in\R^n$. A direct calculation yields the
following formula:
 \bdm
 C_m\ =\ \frac{1}{m!}\, \left( \frac{k^2}{4}+\frac{k}{4}\sqrt{k^2+8m}\right)^m
 \exp \left( \frac{k^2}{4}+\frac{k}{4}\sqrt{k^2+8m} - 
 \frac{1}{k^2} \left(\frac{k^2}{4}+\frac{k}{4}\sqrt{k^2+8m}\right)^2\right)\, .
 \edm
We are only interested in the asymptotics of $C_m$. We will thus ignore all
constant factors not depending on $m$. In this sense, we obtain
 \bdm
 C_m\ \approx\ \frac{1}{m!}\, 
 \left( \frac{k^2}{4}+\frac{k}{4}\sqrt{k^2+8m}\right)^m
 \exp \left(\frac{k}{8}\sqrt{k^2+8m} - \frac{k^2+8m}{16}\right)\, .
 \edm
The Stirling formula $m!\approx \sqrt{m}\,m^me^{-m}$ allows us to rewrite
the asymptotics of $C_m$:
 \bdm
 C_m\ \approx\ \frac{1}{\sqrt{m}\,m^m}
 \left( \frac{k^2}{4}+\frac{k}{4}\sqrt{k^2+8m}\right)^m
 \exp \left(\frac{k}{8}\sqrt{k^2+8m} +\frac{m}{2}\right)\, .
 \edm
Since
 \bdm
 \lim_{m\ra\infty} (\sqrt{k^2+8m} -\sqrt{8m} )\ =\ 0\, ,
 \edm
we can furthermore replace $\sqrt{k^2+8m}$ by  $2\sqrt{2m}$:
 \bdm
  C_m\ \approx\ \frac{1}{\sqrt{m}\,m^m}
 \left( \frac{k^2}{4}+\frac{k}{4}\sqrt{k^2+8m}\right)^m
  \exp \left(\frac{k}{4}\sqrt{2m} +\frac{m}{2}\right)
 \ =:\ e^{C^*_m}
 \edm
with
 \bdm
 C^*_m\ =\ m\ln \left( \frac{k^2}{4}+\frac{k}{4}\sqrt{k^2+8m}\right) 
 +\frac{k}{2\sqrt{2}}\sqrt{m}+\frac{m}{2} -m\ln (m) -\frac{1}{2}\ln (m)\, .
 \edm
In case $m$ is sufficiently large with respect to $k$, we can estimate 
$\ln (k^2/4+k/4\cdot\sqrt{k^2+8m})$
by $\frac{1}{2}\ln (m) +\alpha$ for some constant $\alpha$:
 \bea[*]
 C^*_m & \lessapprox & \frac{m}{2}\ln(m)+\alpha m +\frac{k}{2\sqrt{2}}\sqrt{m}+
 \frac{m}{2} -m\ln (m) -\frac{1}{2}\ln (m)\\
  & \leq & - \frac{m}{2}\ln(m) +(\alpha+1/2)m +\frac{k}{2\sqrt{2}}\sqrt{m}\\
  & \leq & - \frac{m}{2}\ln(m) +(\alpha+1/2+\frac{k}{2\sqrt{2}})m\\
  & =    & m\left(\alpha + 1/2+\frac{k}{2\sqrt{2}}-\frac{1}{2}\ln(m)\right)\, .
 \eea[*]
Finally, $C_m=\exp (C^*_m)$ converges to zero.\qed

\begin{prop}\label{infty-dicht}
Denote by $\Pol[\R^n]$ the ring of all polynomials on $\R^n$. Then the
linear space
$\Sigma_{\infty}\ :=\ \Pol[\R^n] \cdot e^{-r^2}$ is dense
in the space $\Contoinfty{S^n}$ of all continuous functions on
$S^n=\R^n\cup\{\infty\}$ vanishing at infinity.
\end{prop}
\proof
Suppose the closure $\overline{\Sigma_{\infty}}$ of the linear space
$\Sigma_{\infty}$ does not coincide with $\Contoinfty{S^n}$. Then the
Hahn-Banach Theorem implies the existence of a linear continuous
functional $L:\Conto{S^n}\ra\R$ such that
\begin{enumerate}
\item $L\restr_{\Sigma_{\infty}}=0$;
\item $L(g_0)\neq 0$ for at least one $g_0\in\Contoinfty{S^n}$.
\end{enumerate}
According to Riesz' Theorem (see \cite[Ch.6, p.129 ff.]{Rudin2}), $L$
may be represented by two regular Borel measures $\mu_+$, $\mu_-$ on $S^n$:
 \bdm
 L(f)\ =\ \int_{S^n}f(x)\,\de \mu_+(x) - \int_{S^n}f(x)\,\de \mu_-(x)\, .
 \edm
In particular, $\mu_+$ and $\mu_-$ are finite.
The first property $L\restr_{\Sigma_{\infty}}=0$ of $L$ implies
 \bdm
 \int_{S^n}e^{-r^2}p(x)\,\de \mu_+(x)\ =\
 \int_{S^n}e^{-r^2}p(x)\,\de \mu_-(x)
 \edm
for any polynomial $p(x)$. Let us introduce the measures
$\nu_{\pm}=e^{-r^2}\mu_{\pm}$ on the subset $\R^n\subset S^n$. Then
 \bdm
 \int_{\R^n}p(x)\,\de \nu_+(x)\ =\ \int_{\R^n}p(x)\,\de \nu_-(x)
 \edm
holds and remains true for any complex-valued polynomial. 
We may thus choose $p(x)=p_m(x)$ as in the previous lemma
 \bdm
 p_m (x)\ =\ 
 \sum_{\alpha=0}^{m-1} i^{\alpha}\lan k,x \ran ^{\alpha}/ \alpha !\, .
 \edm
But, then
 \bdm
 \int_{S^n}p_m(x)e^{-r^2}\,\de \mu_+(x)\, =\, 
 \int_{\R^n}p_m(x)\,\de \nu_+(x)\, =\, \int_{\R^n}p_m(x)\,\de \nu_-(x)\, =\,
 \int_{S^n}p_m(x)e^{-r^2}\,\de \mu_-(x)
 \edm
together with the uniform convergence of  $p_m(x)e^{-r^2}$
to $e^{i\lan k,x\ran}e^{-r^2}$ implies
 \bdm
 \int_{S^n}e^{i\lan k,x\ran}e^{-r^2}\,\de \mu_+(x)\ =\
 \int_{S^n}e^{i\lan k,x\ran}e^{-r^2}\,\de \mu_-(x)\, ,
 \edm
i.e.,
 \bdm
 \int_{\R^n}e^{i\lan k,x\ran}\,\de \nu_+(x)\ =\
 \int_{\R^n}e^{i\lan k,x\ran}\,\de \nu_-(x)\, .
 \edm
Therefore, the Fourier transforms of the measures $\nu_+$ and $\nu_-$ 
coincide. Consequently, by Bochner's Theorem (see 
\cite[Ch.XIX, p.774 ff.]{Maurin2}) we conclude that $\nu_+=\nu_-$ on $\R^n$. 
The linear functional $L:\Conto{S^n}\ra\R$ must thus be the evaluation of a
function at infinity:
 \bdm
 L(f)\ =\ c\cdot f(\infty)\, ,
 \edm
a contradiction to the existence of a function
$g_0\in\Contoinfty{S^n}$ satisfying $L(g_0)\neq 0$.
\qed
%
\section{The main result}
\begin{thm}
Let the closed subset $M\subset\R^n$ be a smooth submanifold satisfying the
polynomial volume growth condition.
Then the ring $\Aff{\R}{M}$ of all polynomials on $M$ is a dense subspace of the
Hilbert space $L^2(M,e^{-r^2}\de\mu)$.
\end{thm}
\proof
Consider the one-point-compactification $\hat{M}\subset S^n$ of $M\subset\R^n$.
Then Proposition \ref{infty-dicht} of Section $3$ implies that
 \bdm
  \Sigma_{\infty}(\hat{M})\ :=\  \Aff{\R}{M} \cdot e^{-r^2/4}
 \edm
is dense in $\Contoinfty{\hat{M}}$.  We introduce the measure
$\de\nu=e^{-r^2/2}\de\mu$, where $\de\mu$ is the volume form of $M$.
Since
 \bdm
  \int_M \de\nu\ =\ \int_M e^{-r^2/2}\de\mu\ =\
 \int_M (e^{r^2/4})^2 e^{-r^2}\de\mu\ =:\ V\ <\ \infty \, ,
 \edm
$\de\nu$ defines a regular Borel measure $\de\hat{\nu}$ on $\hat{M}$
(by setting $\de\hat{\nu}(\infty)=0$). Therefore, the algebra
$\Contoinfty{\hat{M}}$ 
of all continuous functions on $\hat{M}$ vanishing at infinity is dense in
$L^2(\hat{M},\de\hat{\nu})$:
 \bdm
 \overline{\Contoinfty{\hat{M}}}\ =\ L^2(\hat{M},\de\hat{\nu})\, .
 \edm
For any function $f$ in $L^2(M,e^{-r^2}\de\mu)$ we have
 \bdm
 \int_M | fe^{-r^2/4}|^2 e^{-r^2/2}\de\mu\ =\
 \int_M | f|^2 e^{-r^2}\de\mu\ <\ \infty
 \edm
and, therefore, $fe^{-r^2/4}$ lies in $L^2(\hat{M},\de\hat{\nu})$. Thus,
for a fixed $\eps>0$, there exists a function $g\in \Contoinfty{\hat{M}}$ 
such that
 \bdm
 \int_M  | fe^{-r^2/4} - g(x)|^2 e^{-r^2/2}\de\mu\ <\ \eps/2\, .
 \edm
According to Proposition \ref{infty-dicht} we can find a polynomial
$p(x)\in\Aff{\R}{M}$ approximating $g$:
 \bdm
 \sup_{x\in\hat{M}} |g(x) - p(x) e^{-r^2/4} |^2 \ <\ \eps/2V\, .
 \edm
Using the inequality $\|x+y\|^2\leq 2\|x\|^2+ 2\|y\|^2$ we conclude
 \bdm
 \int_M |f(x)e^{-r^2/4} - p(x)e^{-r^2/4} |^2 e^{-r^2/2}\de\mu\ <\ \eps\, ;
 \edm
but this is equivalent to
 \bdm
 \int_M | f(x) - p(x) |^2 e^{-r^2}\de\mu\ <\ \eps\, .
 \edm
\qed
\section{Examples and final remarks}
We shall give a few simple examples. Notice that we recover, of course, 
that the polynomials are dense in $L^2(\R^n,e^{-r^2}\de\mu)$
(Hermite polynomials) or in $L^2(M,\de\mu)$ for any compact submanifold 
(Legendre polynomials in case $M=[-1,1])$.
\begin{exa}
Consider a revolution surface in $\R^3$ defined by two
polynomials $f,\, h$ 
 \bdm
 \left\{\ba{l}x=f(u_1)\cos u_2\\ y=f(u_1) \sin u_2 \\ z= h(u_1)\ea\right. ,
 \quad f(u_1) > 0\, ,\quad (u_1,u_2)\in \R\x [0,2\pi]\, .
 \edm
Then we have $\de\mu=f\sqrt{{f'}^2+{h'}^2}\, \de u_1\de u_2$ and $r^2=f^2+h^2$,
and thus obtain
 \bdm
 \Aff{\R}{f \cos u_2,\, f \sin u_2,\, h}\ 
 \text{ is dense in }\
 L^2({\R\x [0,2\pi]},e^{-(f^2+h^2)}f\sqrt{{f'}^2+{h'}^2}\, \de u_1\de u_2) \, .
 \edm
In the special case of a cylinder, i.e.\ $f=1$, $h=u_1$, this reduces to the
well known fact that the ring
 \bdm
 \Aff{\R}{u_1,\, \cos u_2,\, \sin u_2,}=\Aff{\R}{u_1}\ox
 \Aff{\R}{\cos u_2,\, \sin u_2}
 \edm
is indeed dense in the Hilbert space
 \bdm
 L^2({\R\x [0,2\pi]},e^{-u_1^2}\,\de u_1\de u_2)=
 L^2(\R,e^{-u_1^2}\,\de u_1)\ox L^2([0,2\pi],\de u_2)\, .
 \edm
%
\end{exa}
\begin{exa}
Let $F:\C\ra\C$ be a polynomial and consider the surface defined by
 \bdm
 f:\ \C\lra \R^3,\quad f(z)=(x,y,|F(z)|),\quad z=x+iy\, .
 \edm
Then one checks that $\de\mu=\sqrt{1+|F'|^2}\, |\de z|^2$ and
$r^2=|z|^2+|F(z)|^2$. Thus the following holds:
 \bdm
 \overline{\Aff{\R}{x,y,|F(z)|}}\ =\ 
 L^2(\R^2,e^{-(|z|^2+|F(z)|^2)}\sqrt{1+|F'|^2}\, |\de z|^2)\, .
 \edm
Let us study the polynomial $F=z^{2k}$ in more detail. Here the
coordinate ring coincides with the usual polynomial ring $\Aff{\R}{x,y}$ 
in two variables, and  thus we have proved that these are dense in
 \bdm
 L^2(\R^2,e^{-(|z|^2+|z|^{4k})}\sqrt{1+4k^2|z|^{2(2k-1)}}\, |\de z|^2)\, .
 \edm
%
\end{exa}
\begin{exa}
We finish with a one-dimensional example: the graph 
$M= \{ (x,f(x)\}$ of a polynomial $f: \R\ra\R^n$. Then
$\de\mu=\sqrt{1+\|f'\|^2}\,\de x$, and we obtain
 \bdm
 \overline{\Aff{\R}{x}}\ =\ 
 L^2(\R,e^{-(x^2+\|f(x)\|^2)}\sqrt{1+ \| f'\|^2}\, \de x)\, .
 \edm
\end{exa}
\begin{rem}
The main result raises an interesting analogous problem in complex analysis
which, to our knowledge, is still open. It is well known that the
polynomials on $\C^n$ are dense in the Fock- or Bergman space 
 \bdm
 \mathcal{F}(\C^n)\ :=\ \{ f\in L^2(\C^n,e^{-r^2}\de\mu)\ |\ 
 f\text{ holomorphic }\}\, .
 \edm
Furthermore, a theorem by Stoll (see \cite{Stoll64a}, \cite{Stoll64b}) states
that from all complex analytic submanifolds $N$ of $\C^n$, those with
polynomial growth are \emph{exactly} the algebraic ones, and thus the
only ones for which the elements of the coordinate ring
are square-integrable with respect to the Gaussian measure.
It is then common to study the space
 \bdm
 \mathcal{F}(N)\ :=\ \{ f\in L^2(N,e^{-r^2}\de\mu)\ |\ 
 f\text{ holomorphic }\}\, ,
 \edm
but we were not able to find any results on whether $\Aff{\C}{N}$
is dense herein.
\end{rem}
\medskip
\noindent
More elaborate applications of the main result to the situation where
$M$ carries a reductive algebraic group action will be discussed
by the authors in some forthcoming works (see e.g.\ \cite{Agricola98}).
In this case, one can decompose the ring $\Aff{\R}{M}$ into
isotypic components and, via Theorem $1$, one obtains a decomposition
of $L^2(M,e^{-r^2}\de\mu)$ analogous to the classical Frobenius reciprocity.

%
\end{document}